\documentstyle{amsppt}
\magnification=1200
\hsize=150truemm
\vsize=224.4truemm
\hoffset=4.8truemm
\voffset=12truemm

\NoRunningHeads

\define\C{{\bold C}}
\let\thm\proclaim
\let\fthm\endproclaim
 


\newcount\tagno
\newcount\secno
\newcount\subsecno
\newcount\stno
\global\subsecno=1
\global\tagno=0
\define\ntag{\global\advance\tagno by 1\tag{\the\tagno}}

\define\sta{\ 
{\the\secno}.\the\stno
\global\advance\stno by 1}

\define\stas{\the\stno
\global\advance\stno by 1}

\define\sect{\global\advance\secno by 1
\global\subsecno=1\global\stno=1\
{\the\secno}. }

\def\nom#1{\edef#1{{\the\secno}.\the\stno}}
\def\inom#1{\edef#1{\the\stno}}
\def\eqnom#1{\edef#1{(\the\tagno)}}

\newcount\refno
\global\refno=0

\def\nextref#1{
      \global\advance\refno by 1
      \xdef#1{\the\refno}}

\def\bref {\ref\global\advance\refno by 1\key{\the\refno}}

\nextref\BR
\nextref\DC
\nextref\DU
\nextref\II
\nextref\NO
\nextref\RU
\nextref\SI
\nextref\YA
\topmatter

\title 
A geometric second main theorem
 \endtitle

\author  Julien Duval and Dinh Tuan Huynh\footnote""{Laboratoire de Math\'ematiques d'Orsay, Univ. Paris-Sud, CNRS, Universit\'e Paris-Saclay, 91405 Orsay, France \newline julien.duval\@math.u-psud.fr\newline dinh-tuan.huynh\@math.u-psud.fr\newline}
\footnote""{keywords : second main theorem, value distribution theory \newline AMS class. : 32H30.\newline}
\endauthor

\abstract\nofrills{\smc Abstract.} We prove a truncated second main theorem in $P^2(\C)$ for entire curves which cluster on an algebraic curve.
\endabstract

\endtopmatter 

\document

\subhead 1. Introduction \endsubhead

\null
Let $f : \C \to P^2(\C)$ be an entire curve (not contained in a line), and $L_i$ be $q$ lines in general position (configuration without triple points). Value distribution theory aims at comparing the growth of $f$ with the number of its impacts on the lines $L_i$. Let $T(f)$ be the Nevanlinna characteristic function and $N_i(f)$ the counting function (with multiplicities). Precisely  
$$T(f,r)=\int_1^r\text{area}(f(D_t))\frac{dt}{t}$$
$$N_i(f,r)=\int_1^r\text{card} (f(D_t)\cap L_i)\frac{dt}{t}$$
 where $D_r$ is the disc centered at $0$ of radius $r$.
Cartan second main theorem reads as follows (outside a set of finite length) 
$$(q-3)T(f,r)\leq \sum N_i(f,r) + o(T(f,r)).$$ 

\noindent
In fact the multiplicities of intersections can even be truncated at level $2$ (see for instance [\RU]).

It is natural to try and prove a similar estimate with multiplicities truncated at level $1$,
meaning that we really count the number of geometric intersections between the entire curve and the lines.
Denote by $N_i^{[1]}(f,r)$ such a truncated counting function.
This is certainly not true if the entire curve is {\it degenerate} (contained in some algebraic curve). Indeed consider a conic $C$, three of its tangents $L_1,L_2,L_3$, and the line $L_4$ joining $L_1\cap C$ to $L_2\cap C$. Take a non constant entire curve $f(\C)$ in $C$ avoiding $L_1,L_2, L_4$. Then $\sum N_i^{[1]}(f,r)=N_3^{[1]}(f,r)=\frac 1 2 T(f,r) +o(T(f,r))$. Hence the question only makes sense for non degenerate entire curves.

\null
Note that in the context of abelian or semi-abelian varieties such a truncated second main theorem has been proved by K. Yamanoi and his collaborators [\YA],[\NO]. 

\null
Here we focus on entire curves whose cluster set is contained in an algebraic curve. Call such entire curves {\it asymptotically degenerate}. There are plenty of them, for instance entire curves which cluster on a line $L$ [\DC]. Indeed observe first that the total space of the line bundle $O(1)$ on $L$ embeds into $P^2(\C)$. So it suffices to construct an entire curve in $O(1)$ which clusters on the zero section. For this take a ramified cover $p : E\to L$ by an elliptic curve. Then $p^*O(1)$ is a positive line bundle on $E$, it inherits a metric of positive curvature. Consider now the universal cover $\pi : \C \to E$. So $(p\circ\pi)^*O(1)$ also has a metric of (uniformly) positive curvature and is holomorphically trivial. Therefore it has a trivializing holomorphic section with exponential decay, giving our asymptotically degenerate entire curve. This construction can be modified to produce entire curves which cluster on any rational or elliptic curve in $P^2(\C)$.

\null
Here is a consequence of our main result.
\thm{Corollary} Let $f : \C \to P^2(\C)$ be an entire curve, non degenerate but asymptotically degenerate. Consider $q$ lines $L_i$ in general position. Then (outside a set of finite logarithmic length) $(q-3)T(f,r)\leq \sum N_i^{[1]}(f,r) + o(T(f,r))$.
\fthm
Finite logarithmic length means of finite measure for $\frac {dr}{r\log r}$. We call such a set an {\it exceptional set}.

This comes from a more general statement best expressed in terms of Nevanlinna currents. We refer to [\BR] for background on these currents. Let $X$ be a projective surface endowed with a hermitian metric and $f : \C \to X$ a non degenerate entire curve. Define the Nevanlinna characteristic function $T(f,r)$ as before. Consider the means of currents of integration $$T_r=\frac 1 {T(f,r)}\int_1^r[f(D_t)]\frac{dt}{t}$$ and their limits as $r\to+\infty$. Then most of these limits are {\it closed} normalized positive currents. Namely this holds true for limits coming from sequences outside an exceptional set. Actually $L(r)=o(T(r))$ outside this set, where $$L(r)= \int_1^r\text{length}(f(\partial D_t))\frac{dt}{t}.$$

Call such a closed limit a {\it Nevanlinna current} associated to $f$. Now given an algebraic curve $C$ in $X$ and a Nevanlinna current $T=\lim T_{r_n}$ we may compute two types of intersection, the algebraic intersection $T.C$ between the two homology classes, and the geometric intersection (without multiplicities) defined by (after extraction) $$i(T,C) =\lim \frac 1 {T(f,r_n)}\int_1^{r_n}\text{card}(f(D_t)\cap C)\frac{dt}{t}.$$  

 By the non degeneracy of $f$ and Poincar\'e-Lelong formula we always have $$i(T,C)\leq T.C. \tag 1$$ In particular the homology class of $T$ is numerically effective. 

\null 
Here is our main result.

\thm {Theorem } Let $f :\C \to X$ be an entire curve, non degenerate but asymptotically degenerate, $C$ an algebraic curve in $X$ with only ordinary double points and $T$ a Nevanlinna current for $f$. Then $T.(K_X+C)\leq i(T,C)$.
\fthm

Here $K_X$ stands for the canonical class of $X$.

Strictly speaking we have to restrict ourselves to Nevanlinna currents $T$ which are projections of similar currents in some blow-up $\pi: \tilde X \to X$. This can always be achieved at the cost of enlarging the exceptional set. Indeed if $\tilde f : \C \to \tilde X$ is the lift of $f$ then the limits $\tilde T$ of $$\tilde T_r=\frac 1 {T(f,r)}\int_1^r[\tilde f(D_t)]\frac{dt}{t}$$ are also closed (non normalized) positive currents if $r\to +\infty$ outside an exceptional set. Therefore we may always consider Nevanlinna currents $T=\lim T_{r_n}$ and $\tilde T=\lim \tilde T_{r_n}$ for the same sequence, in which case $\pi_*\tilde T=T$.
The blow-up is determined by the algebraic data. Namely let $D$ be the algebraic curve containing the cluster set of $f$. By adding extra components we may suppose that $D$ also contains $C$. Then $\tilde X$ is obtained by desingularizing $D$. We refer to [\II] for background in birational geometry.

\null 
In case where $X=P^2(\C)$ and $C$ is a union of $q$ lines in general position we get $q-3 \leq i(T,C)$ for any Nevanlinna current coming from a sequence outside an exceptional set. But this is a reformulation of the corollary. 

\null
The proof of the theorem goes first by reducing the statement to the case of $\tilde X$ and then by using Ahlfors theory of covering surfaces in the spirit of [\DU]. Let us enter the details.

\null
\null
\subhead 2. Proof of the theorem \endsubhead

\null\noindent
{\bf Reduction to $\tilde X$ }. The theorem follows from a similar statement in $\tilde X$
$$\tilde T.(K_{\tilde X}+\pi^{-1}(C))\leq i(\tilde T, \pi^{-1}(C)).  \tag 2$$
Here $\pi^{-1}(C)$ stands for the geometric preimage of $C$, meaning the union of the strict transform of $C$ and the exceptional divisors over $C$. On the other side when we speak of $\pi^*(C)$ we take into account the multiplicities of the exceptional divisors. Note first that $i(\tilde T, \pi^{-1}(C))=i(T,C)$. So it is enough to verify that $$  T.(K_X+C)\leq\tilde T.(K_{\tilde X}+\pi^{-1}(C)) .$$ But this is a consequence of the following inequality between classes (meaning that the difference is effective) $$\pi^*(K_X+C) \leq K_{\tilde X}+\pi^{-1}(C)  ,\tag 3$$ as $\tilde T$ is numerically effective and $\pi_*(\tilde T)=T$. Let us check $(3)$. The map $\pi : \tilde X \to X$ is a composition of blow-ups $\pi_i :X_{i+1}\to X_i$. Denote by $C_i$ the geometric preimage of $C$ in $X_i$.  Its singularities are ordinary double points. We have $\pi_i^*C_i\leq \pi_i^{-1}(C_i)+E_i$ where $E_i$ is the new exceptional divisor. This comes from the fact that the multiplicity of $C_i$ at the point we blow up is at most 2. We also have $K_{X_{i+1}}=\pi_i^*K_{X_i}+E_i$. So we end up with $\pi_i^*(K_{X_i}+C_i)\leq  K_{X_{i+1}}+ \pi_i^{-1}(C_i)$ which gives $(3)$ by iterating.

\null
Now $(2)$ follows from the same statement where we enlarge $C$ to $D$
$$\tilde T.(K_{\tilde X}+\pi^{-1}(D))\leq i(\tilde T, \pi^{-1}(D)).  \tag 4$$
Indeed write $\pi^{-1}(D)=\pi^{-1}(C)+F$ where $F$ is effective. We have $$i(\tilde T,\pi^{-1}(C)+F)\leq i(\tilde T,\pi^{-1}(C))+i(\tilde T,F)\leq i(\tilde T,\pi^{-1}(C))+\tilde T.F$$ by $(1)$. So we get $(2)$ out of $(4)$.

\null
We translate now the left hand side of $(4)$ in more geometric terms. Remark that the Nevanlinna current $\tilde T$ is supported in $\pi^{-1}(D)$ as the cluster set of $f$ is contained in $D$. So we have $\tilde T= \sum \nu(d)  [d]$ where the sum runs over the irreducible components of $\pi^{-1}(D)$ and the weights $\nu$ are non negative. Denote by $d^*$ the component $d$ punctured at the intersection points with the other components and $\chi(d^*)$ its (geometric) Euler characteristic. Recall that the components $d$ are smooth and their intersections consist in ordinary double points as we have desingularized $D$ in $\tilde X$. So $-\chi(d^*)=d.(K_{\tilde X}+\pi^{-1}(D))$ by the genus formula. Summing up we get 
$$ \tilde T.(K_{\tilde X}+\pi^{-1}(D))=-\chi(\tilde T)$$
where $\chi(\tilde T)=\sum \nu(d)\chi(d^*)$ is the Euler characteristic of $\tilde T$.  
So it remains to prove
$$-\chi(\tilde T)\leq i(\tilde T, \pi^{-1}(D)).  \tag 5$$

\null\noindent
{\bf Use of Ahlfors theory }. For this we rely on the technique of [\DU]. 

\null
First note that $\tilde T$ is a particular instance of (non normalized) {\it Ahlfors current} in the sense of [\DU]. This means that $\tilde T=\lim \frac {[\Delta_n]}{t_n}$ where $\Delta_n$ is a union of (possibly repeated) holomorphic discs, $t_n$ is comparable to $\text{area}(\Delta_n)$ and $\text{length}(\partial \Delta_n)=o(t_n)$. Indeed it suffices to discretize the logarithmic means involved in the definition of $\tilde T$ by Riemann sums with rational coefficients. So by [\DU] we already know that the components $d$ charged by $\tilde T$ are rational or elliptic. 

\null
Moreover the technique developped in [\DU] tells us that $\Delta_n$ asymptotically behaves on $d$ as a covering of   degree $t_n \nu(d)$ (at least on graphs). Fix from now on $\epsilon>0$. This means that given a graph $\gamma$ in $d$ we are able to lift it to a (nearby) graph $\gamma_n$ in $\Delta_n$ so that $\vert\chi(\gamma_n)- t_n\nu(d)\chi(\gamma)\vert\leq \epsilon t_n$. We restrict ourselves to the components $d$ with $\chi(d^*)\leq0$. Define $\chi^-(\tilde T)$ to be the part of the sum giving $\chi(\tilde T)$ restricted to these $d$ (and $\chi^+(\tilde T)$ the other part). Consider $\gamma$ a bouquet of circles which is a retract of $d^*$ (in case $d^*$ is really punctured). Then $\chi(\gamma)=\chi(d^*)$ and $d^*\setminus \gamma$ consists in discs with one puncture. Putting these graphs together we get $\Gamma$ which can be lifted to $\Gamma_n$ in $\Delta_n$ with $\vert \chi(\Gamma_n)-t_n \chi^-(\tilde T)\vert\leq c\epsilon t_n$. Here $c$ is the number of irreducible components in $\pi^{-1}(D)$. Now recall the following topological fact. If $\gamma$ is a graph in a disc $\Delta$ (or a union of discs), then $\Delta\setminus \gamma$ contains at least $-\chi(\gamma)$ distinct discs.
So $\Gamma_n$ cuts many holomorphic discs $\delta_n$ out in $\Delta_n$. 

\null
The point is that most of these discs intersect $\pi^{-1}(D)$. Let us make this precise.
Note first that most of the discs $\delta_n$ are of area bounded by a sufficiently big constant. Namely $\text{card}\{\delta_n/\text{area}(\delta_n)>A\}\leq \frac 1 A \text{area}(\Delta_n)\leq \epsilon t_n$ for $A$ big enough. From now on we restrict ourselves to the discs $\delta_n$ with $\text{area}(\delta_n)\leq A$. We analyze their limits. By a classical result (see for instance [\SI]) such a limit is either a holomorphic disc ({\it good limit}) or a union of a holomorphic disc and finitely many rational curves ({\it bubbling limit}). Moreover according to [\DU] we have a good control of the discs $\delta_n$ near their boundaries. We know that $\partial \delta_n$ converges toward an essential loop of (a small perturbation of) the bouquet $\gamma$ in some component $d$. So in either case the limit disc $\delta$ has to cover one of the components of $d\setminus \gamma$. We infer firstly that $\text{area}(\delta_n)$ is bounded away from $0$. Hence the number of the discs $\delta_n$ is bounded from above by $O(t_n)$. Secondly $\delta$ intersects another component of $\pi^{-1}(D)$ as all the components of $d\setminus \gamma$ do. In the good case this still holds before the limit for $\delta_n$ by Hurwitz theorem. Indeed in this case the convergence of $\delta_n$ toward $\delta$ is good (it is actually uniform in some parametrization). We get intersections with $\pi^{-1}(D)$. 

\null
To proceed further we switch back to currents.
Denote by $\Cal P$ the compact set of the positive currents $V$ in $\tilde X$ of mass $\leq A$. Consider the current $U_n=\frac 1 {t_n} \sum[\delta_n]$ where the sum runs over the discs $\delta_n$ defined above. It can also be written as $\int_{\Cal P} V d\mu_n(V)$ where $\mu_n$ is a discrete positive measure on $\Cal P$ of bounded mass. We may suppose that $\mu_n$ converges toward a positive measure $\mu$. Our discs counting above translates in  $$\text{mass}(\mu)\geq -\chi^-(\tilde T)-(c+1)\epsilon.$$ Note that the limit current $U=\int_{\Cal P} V d\mu(V)$ is dominated by $\tilde T$. So the support of $\mu$ consists in currents supported in $\pi^{-1}(D)$. By construction it is also built out of limits of sequences of discs $\delta_n$. Denote by $\mu_g$ and $\mu_b$ the restrictions of $\mu$ to the good limits and to the bubbling limits respectively. Our intersection argument above translates in
$$\text{mass}(\mu_g)\leq i(\tilde T, \pi^{-1}(D)).$$ Let us now look at the bubbling limits. Put $U_b=\int_{\Cal P} V d\mu_b(V)$. We still have $U_b\leq \tilde T$. So all the bubbles in the support of $\mu_b$ are components of $\pi^{-1}(D)$. Moreover such a component $d$ cannot contain a bouquet $\gamma$. Indeed if it did we would have been able to lift $\gamma$ in $\delta_n$ before the limit meaning that $\delta_n$ would intersect $\Gamma_n$. Hence $\chi(d^*)>0$. Specializing the inequality $U_b\leq \tilde T$ on each such $d$ we get $$\text{mass}(\mu_b)\leq \chi^+(\tilde T). $$   
Summing up we have $-\chi(\tilde T)-(c+1)\epsilon \leq i(\tilde T, \pi^{-1}(D))$ and $(5)$ follows by letting $\epsilon\to 0$.
\Refs
\widestnumber\no{99} 
\refno=0 
\bref \by M. Brunella \paper Courbes enti\`eres et feuilletages holomorphes \jour
Ens. Math. \vol45\yr1999\pages195--216
\endref
\bref \by B. Freitas Paulo da Costa \paper Deux exemples sur la dimension moyenne d'un espace de courbes de Brody \jour Ann. Inst. Fourier \vol63\yr2013\pages2223--2237
\endref
\bref \by J. Duval \paper Singularit\'es des courants d'Ahlfors \jour Ann. Sci. ENS \vol39\yr2006\pages527--533
\endref
\bref \by S. Iitaka \book Algebraic geometry. An introduction to birational geometry of algebraic varieties \publ Springer \yr 1982 \publaddr Berlin
\endref 
\bref \by J. Noguchi, J. Winkelmann and K. Yamanoi \paper  The second main theorem for holomorphic curves into semi-abelian varieties. II.  \jour Forum Math. \vol20\yr2008\pages469--503
\endref
\bref \by M. Ru \book Nevanlinna theory and its relation to Diophantine approximation  \publ World Scientific  \yr 2001 \publaddr River Edge
\endref 
\bref \by J.-C. Sikorav \book Some properties of holomorphic curves in almost complex manifolds, {\rm in} Holomorphic curves in symplectic geometry
\pages165--189 \publ Birkh\"auser \yr1994 \publaddr Basel
\endref
\bref \by K. Yamanoi \paper  Holomorphic curves in abelian varieties and intersections with higher codimensional subvarieties \jour Forum Math. \vol16\yr2004\pages749--788
\endref 
 
\endRefs

\enddocument